\documentclass[12pt]{amsart}
\usepackage{amscd,amssymb}
\usepackage[graph,frame,poly,arc]{xy}  
\usepackage[plainpages,backref,urlcolor=blue]{hyperref}
\usepackage{color}

\theoremstyle{plain}
\newtheorem{thm}[subsection]{Theorem}

\newtheorem{prop}[subsection]{Proposition}
\newtheorem{cor}[subsection]{Corollary}

\theoremstyle{definition}
\newtheorem{rk}[subsection]{Remark}
\newtheorem{definition}[subsection]{Definition}
\newtheorem{ex}[subsection]{Example}
\newtheorem{conj}[subsection]{Conjecture}

\numberwithin{equation}{section}
\newcommand{\A}{{\mathcal A}}

\newcommand{\al}{{\alpha}}

\newcommand{\Z}{\mathbb{Z}}

\newcommand{\C}{\mathbb{C}}
\newcommand{\PP}{\mathbb{P}}

\newcommand{\N}{\mathbb{N}}

\DeclareMathOperator{\dd}{d}

\begin{document}

\title [Milnor fiber monodromy]
{Computing Milnor fiber monodromy for some projective hypersurfaces}

\author[Alexandru Dimca]{Alexandru Dimca$^1$}
\address{Universit\'e C\^ ote d'Azur, CNRS, LJAD, France }
\email{dimca@unice.fr}

\author[Gabriel Sticlaru]{Gabriel Sticlaru}
\address{Faculty of Mathematics and Informatics,
Ovidius University,
Bd. Mamaia 124, 
900527 Constanta,
Romania}
\email{gabrielsticlaru@yahoo.com }
\thanks{$^1$ This work has been supported by the French government, through the $\rm UCA^{\rm JEDI}$ Investments in the Future project managed by the National Research Agency (ANR) with the reference number ANR-15-IDEX-01.}



\begin{abstract} We describe an algorithm computing the monodromy and the pole order filtration on the top Milnor fiber cohomology of hypersurfaces in $\PP^n$ whose pole order spectral sequence degenerates at the second page. In the case of
hyperplane arrangements and free, locally quasi-homogeneous hypersurfaces, and assuming a key conjecture, this algorithm is much faster than for a  hypersurface as above. Our conjecture is supported by the results due to  L. Narv\' ez Macarro  and  M. Saito on the roots of  Bernstein-Sato polynomials of such hypersurfaces, by
all the examples computed so far, and by one partial result. 
For hyperplane arrangements coming from reflection groups, a surprising symmetry of their pole order spectra on top cohomology is displayed in our examples. We also improve our previous results in the case of plane curves.

\end{abstract}
 
\maketitle


\section{Introduction} \label{sec:intro}

Let $V:f=0$ be a reduced hypersurface  in the complex projective space $\PP^{n}$, defined by a homogeneous polynomial $f \in S=\C[x_0,...,x_n]$, of degree $d$.
Consider the corresponding complement $M=\PP^{n}\setminus  V$, and the global Milnor fiber $F$ defined by $f(x_0,...,x_n)=1$ in $\C^{n+1}$, with monodromy action $h:F \to F$, $h(x)=\exp(2\pi i/d)\cdot x$. A special case of great interest is when $f$ is a product of linear forms, and then $V$ is a hyperplane arrangement $\A$, and the corresponding complement is traditionally denoted by $M(\A)$.
A lot of efforts were made, in the case of hyperplane arrangements most of the time, to determine  the eigenvalues of the monodromy operators
\begin{equation} 
\label{mono1}
h^m: H^m(F,\C) \to H^m(F,\C)
\end{equation} 
with  $1 \leq m \leq n$, see for instance  \cite{A2, PB, BSett, BY, BDS, Cal1, Cal2, CS, DHA, DL3, MP, MPP, PS, Se2,  S2}. However, in most of these papers, either only the monodromy action on $H^1(F, \C)$ is considered, or the results are just sufficient conditions for the vanishing of some eigenspaces 
$H^m(F,\C)_{\lambda}$. These conditions are usually not necessary, see Example \ref{exNF} below.
For complexified real arrangements, an approach to compute the monodromy operators using the associated Salvetti complex is explained in \cite{Cal1, Cal2, SalSet,  Se2}. However, note that the Milnor fibers considered in \cite{Cal1, Cal2} are not the same as the Milnor fibers in our note, but they correspond to the discriminants of some reflection groups.

In this paper we explain an approach  working  for {\it some} hypersurfaces, namely in technical terms for hypersurfaces $V:f=0$ whose pole order spectral sequence $E_*(f)$ described below degenerates at the $E_2$-term. For hyperplane arrangements and free locally quasi-homogeneous hypersurfaces, modulo a basic conjecture that is {\it one of the main contribution of this paper}, see Conjecture \ref{conj2} below, this algorithm is quite efficient. This conjecture is suggested by the fact that,
for these latter hypersurfaces,  the roots of their Bernstein-Sato polynomials enjoy special properties, as proved by  L. Narv\' ez Macarro \cite{NM} and  M. Saito \cite{Sa0}.
In fact our results are either  {\it conjectural}, depending on whether the Conjecture \ref{conj2}
holds (as for instance in Examples \ref{exD4}, \ref{exG}, \ref{exA214} and \ref{exNF}, and for $k$ a resonant value), or {\it certain}, but based on additional information coming from other viewpoints (as for instance in the previous examples, but for $k$ non-resonant, see Remark \ref{rkprim1} on this point, or in Examples \ref{exFA}, \ref{exBraid}, \ref{exD3}, \ref{exGEN} and \ref{exGEN2}).

Our computation gives not only the dimensions of the eigenspaces $H^n(F,\C)_{\lambda}$ of the monodromy, but also the dimensions of the graded
pieces $Gr^p_PH^n(F,\C)_{\lambda}$, where $P$ denotes the pole order filtration on the  cohomology group $H^m(F,\C)$, see section 2 below for the definition. The dimensions of the eigenspaces $H^m(F,\C)_{\lambda}$ for $m<n$ can then be computed by decreasing induction on $n$, using a generic linear section and the formula \eqref{Euler} below.

In the case $n=2$, this approach was already described in \cite{DStFor, DStMFgen} in the case of a reduced plane curve $C:f=0$. However, even in this case, we bring here valuable new information, see Proposition \ref{propFP} and Corollary \ref{corcurve}. These two results have short and rather direct proofs, but their consequences for the practical computations are substantial, and hence we regard them as {\it main results} of our paper.

Assume now that $n>2$ and let $H \subset \C^{n+1}$ be a generic hyperplane with respect to the hypersurface $V$, passing through the origin. Let $V_H=V \cap H$ be the corresponding hyperplane section of $V$ in $\PP(H)=\PP^{n-1}$, and denote by $F_H$ the corresponding Milnor fiber in $H=\C^n$ and by 
$$h_H^m: H^m(F_H,\C) \to H^m(F_H,\C)$$ 
the associated monodromy operators. Then it is known that the obvious inclusion $\iota_H: F_H \to F$ induces isomorphisms $H^m(F,\C) = H^m(F_H,\C)$ for $m=1,2,..., n-2$, as well as a monomorphism
\begin{equation} 
\label{mono}
\iota_H^*:H^{n-1}(F,\C) \to H^{n-1}(F_H,\C),
\end{equation} 
see for instance \cite{D1}, which are compatible with the monodromy operators.
Consider  the Alexander polynomials of $V$, which are just the characteristic polynomials of the monodromy, namely
\begin{equation} 
\label{Delta}
\Delta^j(V)(t)=\det (t\cdot Id -h^j|H^j(F,\C)),
\end{equation} 
for $j=0,1,...,n$, denoted by $\Delta^j(\A)(t)$ in the case $V=\A$. It is clear that  one has $\Delta^0(V)(t)=t-1$, and moreover
\begin{equation} 
\label{Euler}
\Delta^0(V)(t)\Delta^1(V)(t)^{-1}\cdots \Delta^n(V)(t)^{(-1)^n}=(t^d-1)^{\chi(M)},
\end{equation} 
where  $\chi(M)$ denotes the Euler characteristic of the complement $M$, see for instance \cite[Proposition 4.1.21]{D1}. When $V$ is a hyperplane arrangement $\A$, the Euler characteristic 
$\chi(M(\A))$ is easily computable from the intersection lattice $L(\A)$, see \cite{DHA, OT}. By induction, assume that we know how to compute the characteristic polynomials
$\Delta^j(V_H)(t)$ for $j=0,1,...,n-1$. It follows that 
\begin{equation} 
\label{Delta2}
\Delta^j(V)(t)=\Delta^j(V_H)(t),
\end{equation} 
for $j=0,1,...,n-2$, and hence, in view of the formula \eqref{Euler}, it is enough to determine the top degree Alexander polynomial $\Delta^n(V)(t)$ and the Euler characteristic $\chi(M)=n+1 -\chi(V)$. {\it The computation of this Alexander polynomial 
$\Delta^n(V)(t)$ is the main aim of this paper}. For the computation of the Euler characteristic  $\chi(V)$, see \cite[Corollary 2]{MSS}.

Here is in short how we proceed. Let $\Omega^j$ denote the graded $S$-module of (polynomial) differential $j$-forms on $\C^{n+1}$, for $0 \leq j \leq n+1$. The complex $K^*_f=(\Omega^*, \dd f \wedge)$ is just the Koszul complex in $S$ of the partial derivatives $f_0, f_1 ... f_n$ of the polynomial $f$ with respect to $x_0, x_1,...,x_n$. The general theory says that there is a spectral sequence $E_*(f)$, whose first term $E_1(f)$ is computable from the cohomology of the Koszul complex $K^*_f$ and whose limit
$E_{\infty}(f)$ gives us the action of the monodromy operator  on the graded pieces of the reduced cohomology $\tilde H^*(F,\C)$ of the Milnor fiber with respect to the pole order filtration $P$, see \cite{Dcomp, DS1, Sa3, Sa4} as well as  \cite[Chapter 6]{D1}. In this note we present  an algorithms to compute the second page of the spectral sequence $E_*(f)$. 
Several examples computed so far suggest the following.
\begin{conj}
\label{conj1}
The spectral sequences  $E_*(f)$  degenerates at the $E_2$-term for any hyperplane arrangement  $\A:f=0$ and any free locally quasi-homogeneous hypersurface $V:f=0$ in $\PP^n$.
\end{conj} 
For the moment it is not clear how to prove this conjecture, not even how to check that it holds in a specific example. For a related property, {\it extremely useful for performing our computations}, see Definition \ref{def1},  Remark \ref{rkprim1} and Conjecture \ref{conj2}. This  property holds in all the cases where we dispose of enough additional information to compute the monodromy operators, see Examples \ref{exFA}, \ref{exBraid}, \ref{exD3}. It also holds for some irreducible non-free surfaces, see Examples \ref{exGEN} and \ref{exGEN2}. Theorem \ref{thmlog1} gives some theoretical support for Conjecture \ref{conj2}, and is a final {\it main result} in our paper.
Conjecture \ref{conj1} can also be regarded as an extension of  the following recent deep result due to M. Saito \cite{Sa3}. 

\begin{thm}
\label{thmconj0}
If a hypersurface $V:f=0$ in $\PP^n$ has only isolated singularities,
then
the spectral sequence $E_*(f)$ degenerates at the $E_2$-term if and only if these  singularities are weighted homogeneous. In particular, Conjecture \ref{conj1} holds for $n=2$.
\end{thm}

From a different point of view, Conjecture \ref{conj1}  can  be regarded as a special case of a general conjecture for singular projective hypersurfaces going back to H. Terao \cite{Terao78}, 
and saying that $E_2(f)=E_{\infty}(f)$ always holds. This conjecture is known to fail in general, e.g. by Theorem \ref{thmconj0} above or by looking at the surface $V':f'=0$ introduced at the end of  Example \ref{exGEN}, see also \cite{Dcomp}. The remarkable fact pointed out in our paper is that Terao's Conjecture seems to hold {\it in a stronger form} for any hyperplane arrangement.

In the final section several examples of plane arrangements in $\PP^3$, 
as well as  examples of  (free, locally quasi-homogeneous  or general) surfaces in $\PP^3$, are considered to illustrate the method.  For reflection arrangements, the pole order spectrum $Sp^0_P(f)$ of the top cohomology group $H^3(F, \C)$ has a surprising symmetry property, see Remark \ref{rkKEY4}, which is not present in the case of other arrangements considered in Examples \ref{exFA} and \ref{exNF}.
There is no explanation for this symmetry for the moment, just a possible analogy to the formula \eqref{eqMN} verified by the Bernstein-Sato polynomial of a free, locally quasi-homogeneous hypersurface.

The computations in this note were made using the computer algebra system  Singular \cite{Sing}.
The corresponding codes are available on request.

\medskip

We thank Morihiko Saito for very useful discussions related to the subject and the presentation of this note, see in particular Remark \ref{rkF=P} and Remark \ref{nonfreearr}.

\section{Gauss-Manin complexes, Koszul complexes, and Milnor fiber cohomology} \label{sec2}

Let $S$ be the polynomial ring $\C[x_0,...,x_n]$ with the usual grading and consider a reduced homogeneous  polynomial $f \in S$ of degree $d$. The graded Gauss-Manin complex $C_f^*$ associated to $f$ is defined by taking $C_f^j=\Omega^j[\partial_t]$, i.e. formal polynomials in $\partial_t$ with coefficients in the space of differential forms $\Omega^j$, where $\deg \partial_t=-d$ and the differential
$\dd: C_f^j \to C_f^{j+1}$ is $\C$-linear and given by 
\begin{equation} 
\label{difC}
 \dd (\omega \partial_t^q)=(\dd \omega)\partial_t^q-(\dd f \wedge \omega) \partial_t^{q+1},
\end{equation} 
see for more details \cite[Section 4]{DS1}. 
The complex $C_f^*$ has a natural increasing filtration $P'_*$ defined by
\begin{equation} 
\label{filC}
 P'_qC^j_f=\oplus_{i \leq q+j}\Omega^j\partial_t^i.
\end{equation} 
If we set $P'^q=P'_{-q}$ in order to get a decreasing filtration, then one has
\begin{equation} 
\label{grC}
 Gr^q_{P'}C^*_f=\sigma_{\geq q}(K^*_f((n+1-q)d)),
\end{equation} 
the truncation of a shifted version of the Koszul complex $K^*_f$.
Moreover, this filtration $P'^q$ yields a decreasing filtration $P'$ on the cohomology groups $H^j(C^*_f)$ and a spectral sequence
\begin{equation} 
\label{spsqC}
 E_1^{q,j-q}(f) \Rightarrow H^j(C^*_f).
\end{equation} 
On the other hand, the reduced cohomology $\tilde H^j(F,\C)$ of the Milnor fiber $F:f(x_0,...,x_n)=1$ associated to $f$ has a pole order decreasing filtration $P$, see \cite[Section 3]{DS1}, such that there is a natural identification for any integers $q$, $j$ and $k \in [1,d]$
\begin{equation} 
\label{filH}
 P'^{q+1}H^{j+1}(C^*_f)_k=P^q\tilde H^j(F,\C)_{\lambda},
\end{equation} 
where $\lambda=\exp (-2 \pi ik/d).$ Moreover, the $E_1$-term of the spectral sequence 
\eqref{spsqC} is completely determined by the morphisms of graded $\C$-vector spaces
\begin{equation} 
\label{diff1}
 \dd ' : H^j(K^*_f) \to H^{j+1}(K^*_f),
\end{equation} 
induced by the exterior differentiation of forms, i.e. $\dd ' :[\omega] \mapsto [\dd (\omega)]$.
 More precisely, this spectral sequence $E_*(f)$ is the direct sum of $d$ spectral sequences $E_*(f)_k$,
 for $k \in [1,d]$, where
\begin{equation} 
\label{newspsq}
E_1^{s,t}(f)_k=H^{s+t+1}(K^*_f)_{td+k} \text{ and } \dd_1:E_1^{s,t}(f)_k \to E_1^{s+1,t}(f)_k,  \  \  \dd _1:[\omega] \mapsto [\dd (\omega)].
\end{equation} 
With this notation, one has
\begin{equation} 
\label{limit0}
E_{\infty}^{s,t}(f)_k=Gr_P^s\tilde H^{s+t}(F,\C)_{\lambda }.
\end{equation} 
Since the Milnor fiber $F$ is a smooth affine variety, its cohomology groups $H^m(F, \C)$ have a decreasing Hodge filtration $F$ coming from the mixed Hodge structure constructed by Deligne, see \cite{PeSt}. The two filtrations $P$ and $F$ are related by the inclusion
\begin{equation} 
\label{PFincl}
 F^sH^{m}(F,\C) \subset P^s  H^m(F,\C),
\end{equation} 
for any integers $s,m$, see formula $(4.4.8)$ in \cite{DS1}. This inclusion and the equality $F^0H^{m}(F,\C)=H^{m}(F,\C)$, imply the vanishing
\begin{equation} 
\label{limit}
E_{\infty}^{s,t}(f)_k=0
\end{equation} 
for any $s<0$ or $t<0$, in other words the limit page of the spectral sequence is contained in the first quadrant. Moreover, for $k=d$, one has in addition $E_{\infty}^{0,n}(f)_d=0$, see \cite[Proposition 5.2]{BS0}.
\begin{ex} 
\label{exsmooth}
If $F$ is the Milnor fiber associated to a smooth hypersurface $V$ in $\PP^n$, then the  inclusion \eqref{PFincl} becomes an equality, see for instance \cite{Steen} and \cite[Example 6.2.13]{ D1}. Moreover one has 
\begin{equation} 
\label{limit2}
E_{\infty}^{s,t}(f)_k=0 \text{ for } s+t \ne n  \text{ and } E_{\infty}^{n-t,t}(f)_k=\mu(td+k-n-1),
\end{equation} 
where $\mu(a)$ is the coefficient of $t^a$ in the polynomial
$$\left( \frac{1-t^{d-1}}{1-t}\right)^{n+1}.$$
In particular, $\mu(a) \ne 0$ if and only if $1 \leq a \leq (n+1)(d-2)$. 
\end{ex} 

One has the following result,
the second part of which answers positively a conjecture made in Remark 2.9 (i) in \cite{DStMFgen}.

\begin{prop}
\label{propFP} Let $V$ be a  hypersurface  in $\PP^n$ and $H$ a generic hyperplane.
Consider the linear inclusion $\iota_H:F_H \to F$ defined in the Introduction. Then the induced morphisms
$\iota_H^*:H^{m}(F,\C) \to H^{m}(F_H,\C)$ are strictly compatible with the Hodge filtration $F^p$ and compatible with the pole order filtration $P^p$.
Moreover, one has the following.
\begin{enumerate}

\item If $V$ has only isolated singularities, then the Hodge filtration $F^p$ and the pole order filtration
$P^p$ coincide on $H^{n-1}(F,\C)$.

\item For any hypersurface $V$, in particular for any hyperplane arrangement $\A$, the Hodge filtration $F^p$ and the pole order filtration
$P^p$ coincide on $H^{1}(F,\C)$.

\end{enumerate}
\end{prop}

\proof
Since $\iota_H:F_H \to F$ is a regular mapping, the strict compatibility of $\iota_H^*$ with the Hodge filtration $F^p$ is well known, see \cite{De}. In particular, for $m=n-1$ and $V$ with isolated singularities, since $\iota_H^*$ is injective, this means that a cohomology class $\al \in H^{n-1}(F,\C)$ satisfies 
\begin{equation} 
\label{C1}
\al \in F^pH^{n-1}(F,\C) \text{ if and only if } \iota_H^*(\al )\in F^pH^{n-1}(F_H,\C).
\end{equation} 
The compatibility of $\iota_H^*$ with the pole order filtration $P^p$ means that
\begin{equation} 
\label{C2}
\al \in P^pH^{m}(F,\C) \text{ implies } \iota_H^*(\al )\in P^pH^{m}(F_H,\C).
\end{equation} 
This property comes from the fact that  $\iota_H^*$ induces a morphism
$\iota_H^*: C^*_f \to C^*_{f_H}$ between the corresponding Gauss-Manin complexes, 
where $f_H$ denotes the restriction of the polynomial $f$ to the hyperplane $H$, thought of here as a hyperplane in $\C^{n+1}$. This morphism preserves the $P'_*$ filtrations introduced in \eqref{filC}, i.e. one clearly has
$\iota_H^* (P'_qC^*_f) \subset P'_qC^*_{f_H},$ for any integer $q$.
 To prove the claim (1), it is enough in view of the inclusion in \eqref{PFincl} to prove the converse inclusion  $P^pH^{n-1}(F,\C) \subset F^pH^{n-1}(F,\C)$ for any $p$.
So take  $\al \in P^pH^{n-1}(F,\C)$. Using \eqref{C2}, it follows that 
$$\iota_H^*(\al )\in P^pH^{n-1}(F_H,\C)=F^pH^{n-1}(F_H,\C).$$
The last equality is due to the fact that the hyperplane section $V_H$ being smooth,  the Hodge filtration $F^p$ and the pole order filtration
$P^p$ coincide on $H^{n-1}(F_H,\C)$, as seen in Example \ref{exsmooth}. We conclude using  \eqref{C1}. The claim (2) follows in a similar way, by taking a generic $(n-2)$-codimensional linear section instead of the hyperplane $H$.

\endproof

\begin{rk}
\label{rkF=P}

(i) We thank Morihiko Saito for teaching us to be strict about the difference between strict compatibility and compatibility of morphisms of filtered objects in the above proof.
This corrects a serious gap in our original presentation of the proof above.
See also \cite[Remark 4.4]{Sa3} for a similar property in a general context.

\noindent (ii) We do not know whether the morphism $\iota_H^*:H^{m}(F,\C) \to H^{m}(F_H,\C)$ is strictly compatible with the pole order filtration $P^p$. Indeed, a morphism of filtered complexes, strictly compatible with the filtrations, does not induce in general a strictly compatible morphism when we pass to cohomology groups, see \cite[Section 1.10]{DS2}

If this strict compatibility holds for $\iota_H^*:H^{m}(F,\C) \to H^{m}(F_H,\C)$,
and if one knows the $P$-filtration on $H^*(F_H,\C)$, in order to determine it on $H^*(F,\C)$, it is enough to determine the $P$-filtration on the top cohomology $H^n(F,\C)$ and to identify the image of $H^{n-1}(F,\C)$ inside $H^{n-1}(F_H,\C)$ under $\iota^*_H$. 

\noindent (iii)
It is known that the equality $F^s=P^s$ does not hold on $H^2(F,\C)$, even when $F$ is the Milnor fiber of a line arrangement in $\PP^2$, see for instance \cite[Remark 2.9.(ii)]{DStMFgen}.
\end{rk}

The following result is a major improvement of Theorem 1.2 in \cite{DStFor}. 

\begin{cor} 
\label{corcurve}
For any  curve $V:f=0$ in $\PP^2$, in order to compute the corresponding Alexander polynomial $\Delta^1(V)$, it is enough to compute  the dimensions $\dim E_2^{1,0}(f)_k$ for any $k=1,2,...,d$. More precisely, let $\lambda= \exp(-2\pi ik/d)$ and  $m(\lambda)$ be the multiplicity of $\lambda$ as a root of the Alexander polynomial $\Delta^1(V)$. Then one has
$$m(\lambda)=\dim E_{2}^{1,0}(f)_k+\dim E_{2}^{1,0}(f)_{d-k},$$
for $1 \leq k<d$ and $m(1)=\dim E_{2}^{1,0}(f)_d$.
\end{cor} 

\proof It is well known, e.g. one can use the proof of Proposition \ref{propFP} above, that $H^1(F,\C)_{\ne 1}$ is a pure Hodge structure of weight 1. For any $\lambda= \exp(-2\pi ik/d) \ne 1$, and with obvious notation, it follows that
$$m(\lambda)=h^{1,0}(H^1(F,\C)_{\lambda})+h^{0,1}(H^1(F,\C)_{\lambda})=h^{1,0}(H^1(F,\C)_{\lambda})+h^{1,0}(H^1(F,\C)_{\overline \lambda}),$$
where ${\overline \lambda}$ denotes the complex conjugate of $\lambda$. On the other hand, we have
$$h^{1,0}(H^1(F,\C)_{\lambda})=\dim Gr_P^1\tilde H^{1}(F,\C)_{\lambda }=\dim E_{\infty}^{1,0}(f)_k=\dim E_{2}^{1,0}(f)_k,$$
where the last equality is obvious. These equalities yield our claim for $\lambda \ne 1$. The claim for
$\lambda=  1$ follows from the fact that $H^{1}(F,\C)_1=H^{1}(M,\C)$ is a pure Hodge structure of weight $(1,1)$.

\endproof 

\begin{rk}
\label{rkSpF}
Consider the $j$-th {\it Hodge spectrum}  of the plane curve $V:f=0$, defined by
\begin{equation} 
\label{sp1}
Sp_F^j(f)=\sum_{\al>0}n^j_{F,f,\al}t^{\al} 
\end{equation} 
for $j=0,1$, where 
$$n^j_{F,f,\al}=\dim Gr_F^pH^{2-j}(F,\C)_{\lambda}$$
with $p=[3-\al]$ and $\lambda=\exp(-2 \pi i \alpha)$. When $V$ is a line arrangement, then simple formulas for the difference
$$Sp_F(f)=Sp_F^0(f)-Sp_F^1(f)$$
are given in \cite{BS}. It follows from Proposition \ref{propFP} and the proof of Corollary \ref{corcurve}, that once we know the dimensions $\dim E_2^{1,0}(f)_k$ for any $k=1,2,...,d-1$, we can compute spectrum $Sp_F^1(f)$, and hence via \cite{BS}, the spectrum $Sp_F^0(f)$ as well.
This gives us precise information on the Hodge structure on $H^2(F,\C)$ in this case.
\end{rk}

\section{Hyperplane arrangements, free locally quasi-homogeneous divisors, and Bernstein-Sato polynomials}
In this section with explain why the limit page of the spectral sequences discussed above
enjoy a very useful property in the case of hyperplane arrangements.
Let $(D,0):g=0$ be a complex analytic hypersurface germ at the origin of $\C^{n+1}$ and denote by
$b_{g,0}(s)$ the corresponding (local) Bernstein-Sato polynomial. If the analytic germ $g$ is given by a homogeneous polynomial, then one can define also the global Bernstein-Sato polynomial $b_g(s)$ of $g$,
and one has an equality $b_g(s)=b_{g,0}(s)$, see for more details \cite{SaSurvey, Sa1, Sa2}.
Let $R_{g,0}$  be the set of roots of the polynomial $b_{g,0}(-s)$. When $g$ is a homogeneous polynomial, we use the simpler notation $R_g=R_{g,0}$.

In this section we consider the case when $g=f$ is the defining equation of a hypersurface  $V$ in $\PP^n$ and denote by $D=CV$  the affine cone over $V$, defined in $\C^{n+1}$ by the equation $f=0$.
Recall  M. Saito's  fundamental results in  \cite[Theorem 2]{Sa1} and \cite[Theorem 1]{Sa0}.

\begin{thm}
\label{thmBS}
Let $V: f=0$ be a  hypersurface in $\PP^n$, let $\alpha >0$ be a rational number 
and set $\lambda=\exp(-2 \pi i \alpha)$.
\begin{enumerate}

\item If $Gr_P^p H^{n}(F,\C)_{\lambda }\ne 0$, where $p=[n+1-\al]$, then $ \al \in R_f$.

\item If the sets $\al+\N$ and  $ \cup_{a \in D, a \ne 0}R_{f,a}$ are disjoints,
then the converse of the assertion $(1)$ holds.
\end{enumerate}

\end{thm}

\begin{thm}
\label{thmBS2}
Let $\A: f=0$ be an arrangement of $d$ hyperplanes in $\PP^n$.
Then 
$$\max R_f <2-\frac{1}{d}.$$
\end{thm}

\begin{rk} 
\label{rkKEY2}
As shown by Narv\' ez Macarro, the Bernstein-Sato polynomial $b_f$ of a free arrangement
$\A:f=0$ satisfies the equality
\begin{equation} 
\label{eqMN}
b_f(s-2)=\pm b_f(-s),
\end{equation} 
see \cite{NM}. This equality  implies that the zero set $R_f \subset (0,2)$ is stable under  the involution 
$\al \mapsto 2-\al$, including the multiplicities of the roots. 
In fact,  the equation \eqref{eqMN} holds for a larger class of free hypersurfaces, namely those of {\it linear Jacobian type}, see \cite{NM}. As noted in \cite[Corollary 4.3]{NM}, for any such free hypersurface one has $R_f \subset (0,2)$. A locally quasi-homogeneous divisor $V: f=0$ in $\PP^n$ is of linear Jacobian type, see \cite[Theorem (1.6)]{NM}, and an example of such a surface in given below, see Example \ref{exD3}. Indeed, it is easy to see that the hypersurface $V$ and its affine cone $D=CV$ (regarded as a germ at the origin), are locally quasi-homogeneous divisors in the same time. Moreover, it is clear that $V$ (as a projective hypersurface, see \cite[Section 8.1]{DHA}) and its affine cone $D=CV$  (regarded as a germ at the origin) are free in the same time.

\end{rk} 

{\it Assume from now on in this section, except in Definition \ref{def1} and in the final subsection \ref{planecurve}, that $V: f=0$ is either a  hyperplane arrangement  in $\PP^n$, or a free locally quasi-homogeneous divisor in $\PP^n$. Let $\delta_{k,d}$ be 1 if $k =d$ and 0 otherwise, and $\lambda= \exp(-2\pi i k/d)$.}

 \begin{cor}
\label{corPfilt}
With the above assumption, one has
$$Gr_P^p H^{n}(F,\C)= 0$$
for any $p \leq n-2$ and $Gr_P^{n-1} H^{n}(F,\C)_1= 0$. In other words,
for any $k=1,...,d$ one has
$$P^{n-1+\delta_{k,d}}H^n(F,\C)_{\lambda }=H^n(F,\C)_{\lambda }.$$
In particular $E^{n-t,t}_{\infty}(f)_k=0$ for  $t>1-\delta_{k,d}$.

\end{cor}

\proof
Assume $Gr_P^p H^{n}(F,\C)_{\eta } \ne 0$ for some $p$ and some $\eta=\exp(-2 \pi i \alpha)$ with  $p \leq n+1 -\al<p+1$.  Then Theorem \ref{thmBS} (1) implies that $\al \in R_f$ and hence  using  Theorem \ref{thmBS2} or Remark \ref{rkKEY2}, we get $2>\al >n-p.$ If $p\leq n-2$, or if $p=n-1$ and $\al$ is an integer, then we get a
 contradiction. 

\endproof

\begin{rk}
\label{rkF=P2}
If the morphism $\iota_H^*:H^{m}(F,\C) \to H^{m}(F_H,\C)$ is strictly compatible with the pole order filtration $P^p$,
then one has in addition the following property: 

\medskip

$(\star)$ for a hyperplane arrangement,  $E^{s,t}_{\infty}(f)_k=0$ for  any $s\geq 0$  and $t\geq 2-\delta_{k,d}$.

\medskip

Indeed, the property $(\star)$ is equivalent to $P^{m-1+\delta_{k,d}}H^m(F,\C)_{\lambda }=H^m(F,\C)_{\lambda }$ for any $m$. It follows that in this situation the limit page  $E^{s,t}_{\infty}(f)_k$  is not only contained in the first quadrant, but in fact the non-zero terms are situated on just two horizontal lines when $k \ne d$, namely the lines $t=0$ and $t=1$ (and, respectively, only on the line $t=0$ when $k=d$).

\end{rk}

If we assume Conjecture \ref{conj1} and proceed by induction, it remains to compute the dimension of the terms
$E^{n-q,q}_{2}(f)_k=0$ for any $q\in \{0,1\}$ and $k=1,2,...,d$. It is this property that makes possible the computations in a reasonable amount of time. Since Conjecture \ref{conj1} is difficult to check in practice, we introduce the following notion, motivated by Corollary \ref{corPfilt}.
\begin{definition}
\label{def1}
Let $k,m$ be positive integers
satisfying $1 \leq k \leq d$ and $1 \leq m \leq n$. We say that the hypersurface $V:f=0$ is $(k,m)$-top-computable if $E^{n-t,t}_{\infty}(f)_k=0$ for $m-\delta_{k,d}<t\leq n-\delta_{k,d}$, and 
$$\sum_{t=0}^{m-\delta_{k,d}}\dim E^{n-t,t}_{2}(f)_k =\dim H^n(F,\C)_{\lambda }.$$

\end{definition}
The vanishing of $E^{0,n}_{\infty}(f)_d=0$, which is essential for this definition, follows from \cite[Proposition 5.2]{BS0}.
As an example, a smooth hypersurface $V:f=0$  in $\PP^n$ is $(k,n)$-top-computable for any $k$, but not $(k,n-1)$-top-computable. For  any hypersurface $V:f=0$ and any $t$, one has
$\dim E^{n-t,t}_{2}(f)_k  \geq \dim E^{n-t,t}_{\infty}(f)_k $. Hence, if $V$ is  $(k,m)$-top-computable then $ E^{n-t,t}_{2}(f)_k  = E^{n-t,t}_{\infty}(f)_k $ for $0 \leq t \leq m-\delta_{k,d}$. In particular, the information on the $P$-filtration on $H^n(F,\C)$ given by the second term $E_2$ is also complete in this case. For a hypersurface $V:f=0$ not covered by Corollary \ref{corPfilt}, the simplest way to check the vanishings $E^{n-t,t}_{\infty}(f)_k=0$ for $t>m-\delta_{k,d}$ is to use the vanishings in \eqref{limit} and to check by a direct computation whether $E^{n-t,t}_{2}(f)_k=0$ for $m-\delta_{k,d}<t \leq n-\delta_{k,d}$.
Using this approach, irreducible, non-free surfaces in $\PP^3$ that are still $(k,1)$-top-computable are displayed in Examples \ref{exGEN} and \ref{exGEN2}.
\begin{rk}
\label{rkprim1}
The conditions in Definition \ref{def1} are easy to check as soon as we know $\dim H^n(F,\C)_{\lambda }$ in the case of an arrangement $\A$ of $d$ hyperplanes. Note that the vanishings necessary for the $(k,1)$-top-computability hold by Corollary \ref{corPfilt}.
Let $d'=d/e$ where $e=G.C.D. (d,k)$. Then $\lambda$ is a $d$-root of unity of order $d'$. If there is a hyperplane $H \in \A$, such that for any dense edge $X \subset H$, the number $n_X$ of hyperplanes in $\A$ containing $X$ is not a multiple of $d'$, then it is known that $H^m(F,\C)_{\lambda }=0$ for $m<n$ and $\dim H^n(F,\C)_{\lambda }=|\chi(M(\A)|$, see \cite[Theorem 6.4.18]{D2}, \cite{Lib}. In this case we say that $k$ is {\it non-resonant with respect to the arrangement} $\A$.
Hence, the answer given by the second page of the spectral sequence is  correct  as soon as we know that $k$ is non-resonant and we have the equality
$$\sum_{q=0}^{1-\delta_{k,d}}\dim E^{n-q,q}_{2}(f)_k =|\chi(M)|.$$
The new information given in such a case by the second page of the spectral sequence concerns the pole order filtration on $H^{n}(F,\C)_{\lambda}$. On the other hand, the fact that this equality holds in all the computed cases gives strong support for Conjecture \ref{conj2} below.

\end{rk}
The following is the main conjecture put forth in our paper.
\begin{conj}
\label{conj2}
For any  arrangement $\A: f=0$ of $d$ hyperplanes in $\PP^n$, and  for any free locally quasi-homogeneous divisor $V:f=0$ of degree $d$ in $\PP^n$, the defining polynomial $f$ is $(k,1)$-top-computable
for any  positive integer $k$
satisfying $1 \leq k \leq d$. 
\end{conj}

One has  the following  partial result, saying that Conjecture \ref{conj2} holds for $k=d$.
\begin{thm}
\label{thmlog1}
With the above assumption on $V: f=0$, one has
$$\dim E^{n,0}_{2}(f)_d = \dim E^{n,0}_{\infty}(f)_d =\dim H^n(F,\C)_{1} =\dim H^n(M,\C).$$

\end{thm}

\proof
Let us denote as above by $D$ the affine cone in $X=\C^{n+1}$ over $V$. Then let us denote by $\Omega^*(*D)$ the de Rham complex of rational differential forms on $X$ with poles of arbitrary orders along the hypersurface $D$.
Consider the subcomplex $\Omega^*(D)$ of differential forms in $\Omega^*(*D)$ with logarithmic poles along the hypersurface  $D$. A form $\omega \in \Omega^p(D)$ satisfies, by definition
$$ f \omega \in \Omega^p \text{ and } f \dd \omega \in \in \Omega^{p+1},$$
or, equivalently $\omega=\eta/f$ with $\eta \in \Omega^p$ and $\dd f \wedge \eta$ divisible by $f$.
It is clear that one has a natural identification
$$\Omega^{n+1}(D)_0=H^{n+1}(K_f^*)_d=E^{n,0}_{1}(f)_d,$$
given by $\eta/f \mapsto \eta$, 
where the grading on the $S$-module $\Omega^p(D)$ is the usual one, i.e. $|\eta'/f|=|\eta'|-d$, for any $\eta' \in \Omega^p$. Next the homogeneous component $\Omega^{n}(D)_0$ can be identified by the same map as above to the direct sum
$$ S_{d-n-1}\cdot \omega_n \oplus Syz^n_d$$
where 
$$\omega_n=\sum_{i=0}^n(-1)^ix_i \dd x_0 \wedge \dd x_1 \wedge ... \wedge \widehat {\dd x_i}  \wedge ... \wedge \dd x_n,$$
and 
$$Syz^n_d=\ker \left( \dd f \wedge : \Omega^{n}_d \to \Omega^{n+1}_{2d} \right)=H^{n}(K_f^*)_d=E^{n-1,0}_{1}(f)_d.$$
If $h \in S_{m}$, a direct computation shows that
$$d\left( \frac{h\cdot \omega_n}{f}\right)=(m+n+1-d) \frac{h\cdot \omega'_n}{f} ,$$
where $\omega'_n=\dd x_0 \wedge \dd x_1 \wedge ... \wedge  \dd x_n.$
In particular the image of the differential $d: \Omega^{n}(D)_0 \to \Omega^{n+1}(D)_0$ coincides with the image of the differential $d_1: E^{n-1,0}_{1}(f)_d \to E^{n,0}_{1}(f)_d$, and all the other differentials $d: \Omega^{n}(D)_q \to \Omega^{n+1}(D)_q$ for $q \ne d$ are surjective.
 It follows that
\begin{equation} 
\label{eqK1}
\dim E^{n,0}_{2}(f)_d=\dim H^{n+1}(\Omega^*(D)_0)=\dim H^{n+1}(\Omega^*(D)).
\end{equation} 
On the other hand, we clearly have $\dim H^n(M,\C)= \dim H^{n+1}(X \setminus D, \C)$, see for instance \cite[Prop. 6.4.1]{D2}.

It remains to show that $\dim H^{n+1}(\Omega^*(D))=\dim H^{n+1}(X \setminus D, \C)$.
When $(D,0)$ is a  free locally quasi-homogeneous divisor, this follows from \cite{CMN}.
When $D$ is a hyperplane arrangement, one has to use \cite[Proposition 6.1]{WY}.

\endproof

\subsection{The arbitrary hypersurface case} \label{planecurve}

The algorithm presented below can be applied for any hypersurface $V:f=0$ in $\PP^n$ to compute the terms $E_2^{n-q,q}(f)_k$ of the second page of the above spectral sequences. However, in the general case $Q=qd+k$ takes values up to $(n+1)d$ and not only $2d-1$, which increases dramatically the computer time. Once the second page is computed, then it is a difficult question to decide whether 
\begin{equation} 
\label{eqHappy}
E_2^{n-q,q}(f)_k=E_{\infty}^{n-q,q}(f)_k, 
\end{equation}
and hence whether we have obtained the correct results.
In some cases, one can proceed as follows. Using \eqref{limit}, we see  that
\begin{equation} 
\label{eqHappy2}
\dim H^{n}(F,\C)_{\lambda} \leq  \sum_{q=0,...,n}\dim E_2^{n-q,q}(f)_k
\end{equation}
where $\lambda=\exp(-2 \pi i k/d)$, and equality holds if and only if  one has
the equality \eqref{eqHappy}
for any $q$.
Two examples of such a computation are given below in Examples \ref{exGEN} and \ref{exGEN2}.

\section{The algorithm}

Consider the graded $S-$submodule $AR(f) \subset S^{n+1}$ of {\it all relations} involving the derivatives of $f$, namely
$$r=(r_0,r_1,...,r_n) \in AR(f)_q$$
if and only if  
\begin{equation} 
\label{syz1}
r_0f_0+r_1f_1+ ... +r_nf_n=0
\end{equation} 
 and the polynomials $r_0,r_1,...,r_n$ are in $S_q$.  
Since $S$ is a noetherian ring, the graded $S$-module $AR(f)$ admits a (minimal) system of generators $r^{(j)}$ of Jacobian syzygies, where $j=1,...,g$. 
Assume that
$$r^{(j)}=(r^{(j)}_0,r^{(j)}_1,...,r^{(j)}_n)$$
for $j=1,...,g$ and let $d_j=\deg r^{(j)}_m$, for any $m\in [0,n]$ with $r^{(j)}_m \ne 0$. Assume moreover that
$$d_1 \leq d_2 \leq ... \leq d_g.$$
Such a system of generators can be determined using the software SINGULAR \cite{Sing} or CoCoA \cite{Co}, 
see Remark \ref{rkKEY1}.
To each syzygy $r=(r_0,r_1,...,r_n) \in AR(f)_q$ we can associate an $n$-differential form in $\Omega^n_{q+n}$ by the formula
\begin{equation} 
\label{form1}
\omega(r)=\sum_{i=0}^n(-1)^ir_i \dd x_0 \wedge \dd x_1 \wedge ... \wedge \widehat {\dd x_i}  \wedge ... \wedge \dd x_n.
\end{equation} 
Then equation \eqref{syz1} is equivalent to $\dd f \wedge \omega(r)=0$ and 
\begin{equation} 
\label{form2}
\dd \omega(r)=(\sum_{i=0}^n(r_i)_i )\dd x_0 \wedge \dd x_1 \wedge ...  \wedge \dd x_n,
\end{equation} 
where $(r_i)_i$ denotes the partial derivative of $r_i$ with respect to $x_i$ for $i=0,1,...,n$. It follows that the dimension of $E_2^{n-k,q}(f)_k$, which is the cokernel of the differential 
$$\dd_1:E_1^{n-1-q,q}(f)_k \to E_1^{n-q,q}(f)_k$$
can be computed as follows. Consider the linear mapping
\begin{equation} 
\label{fi4}
\phi_Q:S_{Q-d_1-n} \times ...\times S_{Q-d_g-n}\times S^{n+1}_{Q-d-n} \to S_{Q-n-1},
\end{equation} 
given by 
$$((A_1,...,A_g),(B_0,..., B_n)) \mapsto \sum_{ 0 \leq i \leq n} (
(\sum_{1\leq j \leq g} A_jr^{(j)}_i)_i+B_if_i),$$
for $Q=qd+k$.  This map can be described in a more compact way by using differential forms as follows.
The formula \eqref{form1} implies that the application $\phi_Q$ is nothing else but the map
$$S_{Q-d_1-n} \times ...\times S_{Q-d_n-n}\times  \Omega^{n}_{Q-d} \to \Omega^{n+1}_Q$$
given by
$$((A_1,...,A_n), \eta) \mapsto \dd (\sum_{j=1,n}A_j\omega(r^{(j)}))+\dd f \wedge \eta.$$
Let $R_Q$ be the rank of this linear mapping, which is computed using the software SINGULAR for instance. Then one clearly has
\begin{equation} 
\label{E2dim}
\dim E_2^{n-q,q}(f)_k=\dim S_{Q-n-1}-R_Q= {Q-1 \choose n}- R_Q,
\end{equation} 
for any $q=0,...,n$.

In the case of a hyperplane arrangement $\A:f=0$, or of a free divisor $V:f=0$ of linear Jacobian type,  we can consider only the values $q \leq 1$ if we assume Conjecture \ref{conj2}, while $k \leq d$ by definition. It follows that it is enough to take $Q \leq 2d$ in this case.
In fact,  for a hyperplane arrangement, it is known that $F^nH^n(F,\C)_1=H^n(F,\C)_1=H^n(M(\A), \C)$, which implies that in fact we need to consider only the values $Q \leq 2d-1$ in such a case.

\begin{rk} 
\label{nonfreearr}
In case one does not like to use the generating system of syzygies produced by a computer software, one can proceed in the following more direct way, already considered by us in \cite{DStMFgen} and by Morihiko Saito in \cite{Sa4}.
Let $V: f=0$ be a degree $d$ hypersurface  in $\PP^n$, and for each $Q=qd+k$ consider the linear map
\begin{equation} 
\label{fi5}
\phi^1_Q: \Omega^n_{Q-d} \to \Omega^{n+1}_{Q}, \  \  \ \eta \mapsto \dd f \wedge \eta.
\end{equation} 
Let $\kappa^1(Q)$ be the dimension of the kernel $K^1(Q)$ of this map. Clearly $\kappa^1(Q)=0$ for $Q<d+n$.
Consider next the map
\begin{equation} 
\label{fi6}
\phi^2_Q: \Omega^n_{Q-d} \times \Omega^n_{Q} \to \Omega^{n+1}_{Q} \times \Omega^{n+1}_{Q+d},
\end{equation} 
given by 
$$(\eta_1,\eta_2) \mapsto (\dd f \wedge \eta_1 + \dd (\eta_2), \dd f \wedge \eta_2).$$
Let $\kappa^2(Q)$ be the dimension of the kernel $K^2(Q)$ of this map. Note that one has
$$K^2(Q) \subset \Omega^n_{Q-d} \times K^1(Q+d),$$
and hence the dimension of the vector space $\phi^2_Q(\Omega^n_{Q-d} \times K^1(Q+d))$
is given by
\begin{equation} 
\label{fi7}
R_Q= \dim (\Omega^n_{Q-d} \times K^1(Q+d))-\kappa^2(Q)=\kappa^1(Q+d)-\kappa^2(Q)+ (n+1){Q-d \choose n}.
\end{equation} 
The dimensions  $\kappa^1(Q+d)$ and $\kappa^2(Q)$ can be computed using the software SINGULAR for instance. Then one clearly has the same formula as above, namely
\begin{equation} 
\label{E2dim2}
\dim E_2^{n-k,q}(f)_k=\dim S_{Q-n-1}-R_Q= {Q-1 \choose n}- R_Q.
\end{equation} 
This approach seems to increase the necessary computing time as well as the necessary computer memory substantially. We thank Morihiko Saito for telling us that the algorithm to compute the second page of the spectral sequence using the system of generators is better not only in the case of a free hypersurface, when we had already used it, but also in the general case.
\end{rk}

\begin{rk} 
\label{rkKEY1}
(i) The command $syz(...)$ in the software SINGULAR does not always give a minimal set of generators for the graded 
$S$-module $AR(f)$. 
For the quartic surface discussed in Example \ref{exGEN}, it lists 6 generators for the order 
$rp$, 7 generators for the order $dp$ and 8 generators for the orders $lp$ and $Dp$.
Here $rp =$ reverse lexicographical ordering,
$dp =$ degree reverse lexicographical ordering, 
$lp =$ lexicographical ordering, and
$Dp =$ degree lexicographical ordering. Note also that in some of these listings, the generators are not given with the degrees in increasing order. To get the minimal set of generators one should use the command $minbase(syz( ... ))$.

\medskip

\noindent (ii) To get information on the complexity of computations in the algorithm,
it would be useful to have an upper bound in terms of the geometry of the hypersurface $V:f=0$ on $g$, the minimal number of generators for the $S$-module $AR(f)$, and a lower bound on the minimal degree $d_1$.

Note that for a smooth hypersurface we have the $g=n(n+1)/2$ linear independent Koszul generators of degree $d_1=d-1$.  The equality $g=n(n+1)/2$ holds also for the plane arrangement in Example \ref{exNF} and the singular surface  in Example \ref{exGEN} below.  On the other hand, note that \cite[formula (1.3) and Example 4.3 (i)]{DStEdin}
imply that for a  hypersurface with just one $A_1$ singularity, one has
  $$ g = \frac{n(n+1)}{2}+1.$$
  And \cite[Theorem 1.4]{DStEdin} implies that for a nodal curve $C$ in $\PP^2$, one has
$g \geq r-1$, where $r$ is the number of irreducible components of $C$.
Lower bounds for $d_1$ are known for hypersurfaces having only isolated singularities, see 
\cite[Theorem 9]{DS2} for the case of weighted homogeneous singularities, and \cite[Theorem 2.4]{DAG} for arbitrary isolated singularities.

\medskip

\noindent (iii) Note that \cite[Theorem 1.5 and Example 4.3 (i)]{DStEdin}
imply that for a  hypersurface with just one $A_1$ singularity, one has $d_1=...d_{g-1}=d-1$ and $d_g=n(d-2).$


\end{rk}

\begin{rk}  \label{freehyper}

If $V:f=0$ is a free hypersurface,  then $g=n$ as the $S$-module $AR(f)$ is free. One can verify that the basis given by SINGULAR is correct using Saito's criterion, i.e. the $(n+1)$-square matrix having as the first row $x_0,...,x_n$, and the $j+1$-st row given by $(r^{(j)}_0,r^{(j)}_1,...,r^{(j)}_n)$ for $j=1,...,n$, should have as determinant a constant, non-zero multiple of $f$, see \cite{OT, Te, Yo}. In particular, in this case $d_1+d_2+...+d_n=d-1$.

It is known that, in the case of a hyperplane arrangement $\A:f=0$, the exponents $d_j$ determine the Betti numbers of the complement $M(\A)$, see for instance \cite{OT},
\begin{equation} 
\label{poincare}
\pi(\A,t):=\sum_{i=0}^n b_i(M(\A))t^i=\prod_{j=1}^n(1+d_jt).
\end{equation}

\end{rk} 

\section{Examples}

In this section we consider plane arrangements in $\PP^3$, except in Examples \ref{exD3}, \ref{exGEN} and \ref{exGEN2} where  irreducible quartic surfaces in $\PP^3$ are considered, and we replace the coordinates $x_0,x_1,x_2,x_3$ by $x,y,z,w$. To state the  results, we consider the pole order spectrum defined by
\begin{equation} 
\label{sp1.5}
Sp^0_P(f)=\sum_{\al>0}n_{P,f,\al}t^{\al} 
\end{equation} 
 where 
$$n_{P,f,\al}=\dim Gr_P^pH^{3}(F,\C)_{\lambda}$$
with $p=[4-\al]$ and $\lambda=\exp(-2 \pi i \alpha)$. In view of Corollary \ref{corPfilt} and assuming Conjecture \ref{conj2}, in the case of a plane arrangement the exponents $\al$ with possibly non-zero coefficients $n_{P,f,\al}$ are of the form 
\begin{equation} 
\label{sp2}
\al= \frac{Q}{d}  \text{ and } n_{P,f,\al}=\dim E_2^{n-q,q}(f)_k,
\end{equation} 
where $Q=qd+k$ as above, with $q=0,1$ and $k=1,...,d$.
Note that one has the equality
$$b_n(F)=\sum_{\al>0}n_{P,f,\al}.$$

\begin{ex}[A family of free arrangements]
\label{exFA}
Consider the arrangement
$$\A(p,q): (x^p+y^p)(z^q+w^q)=0.$$
This arrangement is free with exponents $(1,p-1,q-1)$ and the monodromy operators 
$h^m: H^m(F,\C) \to H^m(F,\C)$ can be easily computed using \cite[Theorem 1.4]{DNag} or \cite{Tapp}.
For all the pairs $(p,q)$ we have tested, i.e. $2 \leq p,q \leq d=p+q \leq 12$, the algorithm described above gives the correct result. In other words, the corresponding arrangements $\A(p,q)$ are $(k,1)$-top-computable for all the integers $k$ with $1 \leq k \leq d$. For instance, for the arrangement $\A(4,8)$ we get the following spectrum
\begin{equation} 
\label{spA48}
Sp^0_P(f)=3t^{\frac{6}{12}} +10t^{\frac{9}{12}}+21t^{\frac{12}{12}}+12t^{\frac{15}{12}}+9t^{\frac{18}{12}}+2t^{\frac{21}{12}}.
\end{equation} 
Note that this spectrum is not symmetric with respect to the monomial $21t=21t^{\frac{12}{12}}$, and an equality similar to that in Corollary \ref{corcurve} does not hold for the multiplicities of the roots of $$\Delta^3(\A(4,8))=\Phi_1^{21}\cdot \Phi_2^{12}\cdot \Phi_4^{12}.$$
Here and in the sequel, $\Phi_j$ denotes the $j$-th cyclotomic polynomial. 
\end{ex}

\begin{ex}[The braid arrangement  $A_4$]
\label{exBraid}
The braid arrangement $A_4$ is defined in $\C^5$ by the equation
$$\prod_{0\leq i<j \leq 4}(x_i-x_j)=0.$$
However, this arrangement is not essential. Using the coordinate change
$x=x_1-x_0$, $y=x_2-x_0$, $z=x_3-x_0$, $w=x_4-x_0$, $u=x_0+x_1+x_2+x_3+x_4$,
we see that the essential version of this arrangement is given in $\C^4$, corresponding to $u=0$, by the equation
$$\A: f=xyzw(x-y)(x-z)(x-w)(y-z)(y-w)(z-w)=0.$$
Regarded as an arrangement in $\PP^3$, this arrangement is known to be free with exponents
$(d_1,d_2,d_3)=(2,3,4)$.
Running the algorithm described above and using the fact that \eqref{eqHappy2} is an equality in this case as implied by  Settepanella's results in \cite[Table 2]{Se2}, we get that $\A$ is $(k,1)$-top-computable for any $k\in [1,10]$. In particular, we have
\begin{equation} 
\label{spA5}
Sp^0_P(f)=t^{\frac{4}{10}} +4t^{\frac{5}{10}}+5t^{\frac{6}{10}}+6t^{\frac{7}{10}}+6t^{\frac{8}{10}}+6t^{\frac{9}{10}}+24t^{\frac{10}{10}}+
\end{equation} 
$$+6t^{\frac{11}{10}}+6t^{\frac{12}{10}}+6t^{\frac{13}{10}}+5t^{\frac{14}{10}}+4t^{\frac{15}{10}}+t^{\frac{16}{10}}.$$
Then the formula for the spectrum clearly implies
the following formula for the Alexander polynomial $\Delta^3(\A)$:
\begin{equation} 
\label{AlexA5.3}
\Delta^3(\A)=\Phi_1^{24}\cdot \Phi_2^{8}\cdot \Phi_5^{6}\cdot \Phi_{10}^{6},
\end{equation} 
which coincides of course with the formula given in \cite[Table 2]{Se2}.
It is known that in this case $\Delta^1(\A)=\Phi_1^{9}$, see \cite{MP} or \cite[Table 2]{Se2}. Using the formula \eqref{Euler} and \eqref{poincare}, we get $\chi(M(\A))= -6$ and it follows that
$$\Delta^2(\A)=\Phi_1^{26}\cdot \Phi_2^{2}.$$
This coincides again with the formula given in \cite[Table 2]{Se2}. Any value $k\ne 5, 10$ is non-resonant with respect to the arrangement $A_4$, so for such a $k$, we can get the above results without using Settepanella's results in \cite[Table 2]{Se2}, as explained in Remark \ref{rkprim1}.
\end{ex}

\begin{ex}[The Coxeter arrangement  $D_4$]
\label{exD4}
The  arrangement $D_4$ is defined in $\C^4$ by the equation
$$\A: f=(x^2-y^2)(x^2-z^2)(x^2-w^2)(y^2-z^2)(y^2-w^2)(z^2-w^2)=0.$$
Regarded as an arrangement in $\PP^3$, this arrangement is known to be free with exponents
$(d_1,d_2,d_3)=(3,3,5)$, see \cite{OT}. Running the algorithm described above {\it and assuming Conjecture \ref{conj2} true}, we get
\begin{equation} 
\label{spD4}
Sp^0_P(f)=t^{\frac{4}{12}} +4t^{\frac{5}{12}}+10t^{\frac{6}{12}}+12t^{\frac{7}{12}}+23t^{\frac{8}{12}}+16t^{\frac{9}{12}}+20t^{\frac{10}{12}}+16t^{\frac{11}{12}}+45t^{\frac{12}{12}}+
\end{equation} 
$$+16t^{\frac{13}{12}}+20t^{\frac{14}{12}}+16t^{\frac{15}{12}}+23t^{\frac{16}{12}}+12t^{\frac{17}{12}}+10t^{\frac{18}{12}}+4t^{\frac{19}{12}}+t^{\frac{20}{12}}.$$
This  formula for the spectrum clearly implies
the following formula for the Alexander polynomial $\Delta^3(\A)$:

\begin{equation} 
\label{AlexD4.3}
\Delta^3(\A)=\Phi_1^{45}\cdot \Phi_2^{20}\cdot \Phi_3^{24}\cdot \Phi_4^{16} \cdot \Phi_{6}^{20}\cdot \Phi_{12}^{16}.
\end{equation} 
It is known that in this case $\Delta^1(\A)=\Phi_1^{11}\Phi_3$, see \cite{MP}. Using the formula \eqref{Euler} and \eqref{poincare}, we get $\chi(M(\A))= -16$ and it follows that
$$\Delta^2(\A)=\Phi_1^{39}\cdot \Phi_2^{4} \cdot \Phi_3^{9} \cdot \Phi_6^{4}.$$
Any value $k \ne 2, 4, 6,12$ is non-resonant with respect to the arrangement $D_4$.
\end{ex}

\begin{ex}[The complex reflection arrangement  $\A(3,3,4)$]
\label{exG}
The  hyperplane arrangement $\A(3,3,4)$ is defined in $\C^4$ by the equation
$$\A: f=(x^3-y^3)(x^3-z^3)(x^3-w^3)(y^3-z^3)(y^3-w^3)(z^3-w^3)=0.$$
Regarded as an arrangement in $\PP^3$, this arrangement is known to be free with exponents
$(d_1,d_2,d_3)=(4,6,7)$, see \cite{OT}. Running the algorithm described above {\it and assuming Conjecture \ref{conj2} true}, we get
\begin{equation} 
\label{spG}
Sp^0_P(f)=t^{\frac{4}{18}} +4t^{\frac{5}{18}}+10t^{\frac{6}{18}}+19t^{\frac{7}{18}}+31t^{\frac{8}{18}}+46t^{\frac{9}{18}}+59t^{\frac{10}{18}}+71t^{\frac{11}{18}}+98t^{\frac{12}{18}}+
\end{equation} 
$$+86t^{\frac{13}{18}}+89t^{\frac{14}{18}}+92t^{\frac{15}{18}}+90t^{\frac{16}{18}}+90t^{\frac{17}{18}}+168t^{\frac{18}{18}}+90t^{\frac{19}{18}}+90t^{\frac{20}{18}}+92t^{\frac{21}{18}}+89t^{\frac{22}{18}}+$$
$$+86t^{\frac{23}{18}}+98t^{\frac{24}{18}}+71t^{\frac{25}{18}}+59t^{\frac{26}{18}}+46t^{\frac{27}{18}}+31t^{\frac{28}{18}}+19t^{\frac{29}{18}}+10t^{\frac{30}{18}}+4t^{\frac{31}{18}}+t^{\frac{32}{18}}.$$
This  formula for the spectrum clearly implies
the following formula for the Alexander polynomial $\Delta^3(\A)$:
\begin{equation} 
\label{AlexG.3}
\Delta^3(\A)=\Phi_1^{168}\cdot \Phi_2^{92}\cdot \Phi_3^{108}\cdot \Phi_6^{92} \cdot \Phi_{9}^{90}\cdot \Phi_{18}^{90}.
\end{equation} 
It is known that in this case $\Delta^1(\A)=\Phi_1^{17}\Phi_3$, see \cite{MP}. In fact, a generic plane section $\A_H$ has 42 triple points and 27 nodes, and the result follows also from \cite{BDS}.
Using the formula \eqref{Euler} and \eqref{poincare}, we get $\chi(M(\A))= -90$ and it follows that
$$\Delta^2(\A)=\Phi_1^{94}\cdot \Phi_2^{2} \cdot \Phi_3^{19} \cdot \Phi_6^{2}.$$
Any value $k \ne 3, 6, 9,18$ is non-resonant with respect to the arrangement $\A(3,3,4)$.
\end{ex}
\begin{ex}[The complex reflection arrangement  $\A(2,1,4)$]
\label{exA214}

The  hyperplane arrangement $\A(2,1,4)$ is defined in $\C^4$ by the equation
$$\A: f=xyzw(x^2-y^2)(x^2-z^2)(x^2-w^2)(y^2-z^2)(y^2-w^2)(z^2-w^2)=0.$$
Regarded as an arrangement in $\PP^3$, this arrangement is known to be free with exponents
$(d_1,d_2,d_3)=(3,5,7)$, see \cite{OT}. Running the algorithm described above {\it and assuming Conjecture \ref{conj2} true}, we get
\begin{equation} 
\label{spA214}
Sp^0_P(f)=t^{\frac{4}{16}} +4t^{\frac{5}{16}}+9t^{\frac{6}{16}}+16t^{\frac{7}{16}}+25t^{\frac{8}{16}}+32t^{\frac{9}{16}}+39t^{\frac{10}{16}}+44t^{\frac{11}{16}}+
\end{equation} 
$$+47t^{\frac{12}{16}}+48t^{\frac{13}{16}}+48t^{\frac{14}{16}}+48t^{\frac{15}{16}}+105t^{\frac{16}{16}}+48t^{\frac{17}{16}}+48t^{\frac{18}{16}}+48t^{\frac{19}{16}}+47t^{\frac{20}{16}}+$$
$$+44t^{\frac{21}{16}}+39t^{\frac{22}{16}}+32t^{\frac{23}{16}}+25t^{\frac{24}{16}}+16t^{\frac{25}{16}}+9t^{\frac{26}{16}}+4t^{\frac{27}{16}}+t^{\frac{28}{16}}.$$
This  formula for the spectrum clearly implies
the following formula for the Alexander polynomial $\Delta^3(\A)$:
\begin{equation} 
\label{AlexA214.3}
\Delta^3(\A)=\Phi_1^{105}\cdot \Phi_2^{50}\cdot \Phi_4^{48}\cdot \Phi_8^{48} \cdot \Phi_{16}^{48}.
\end{equation} 
It is known that in this case $\Delta^1(\A)=\Phi_1^{15}$, see \cite{MP}. 
Using the formula \eqref{Euler} and \eqref{poincare}, we get $\chi(M(\A))= -48$ and it follows that
$$\Delta^2(\A)=\Phi_1^{73}\cdot \Phi_2^{2}.$$
Any value $k \ne 6$ is non-resonant with respect to the arrangement $\A(2,1,4)$.
\end{ex}

\begin{rk} 
\label{rkKEY4}
Note that all the spectra coming from reflection groups in Examples \ref{exBraid}, \ref{exD4}, \ref{exG}, \ref{exA214} above enjoy a perfect symmetry with respect the monomial containing $t$, i.e. the coefficients of $t^{\al}$ and $t^{2-\al}$ coincide for all $0<\al <1$. This symmetry might be related to the symmetry of the Bernstein-Sato polynomial $b_f$ of a free arrangement
$\A:f=0$ recalled in \eqref{eqMN}.
 However, note that for some free arrangements as in Example \ref{exFA},  the pole order spectra are not
symmetric, but such arrangements seem to be quite exceptional. Indeed, most of the free arrangements we have tested so far enjoy the above spectrum symmetry property.

\end{rk} 

\begin{ex}[A non free arrangement]
\label{exNF}
Consider the  arrangement $\A$ defined in $\C^4$ by the equation
$$\A: f=xyzw(x+y+z)(y-z+w)=0.$$
This arrangement is far from being free, the  $S$-module $AR(f)$ has 6 generators of degrees $2,2,3,3,3,3$ respectively.
Running the algorithm described above {\it and assuming Conjecture \ref{conj2} true}, we get
\begin{equation} 
\label{spGO}
Sp^0_P(f)=t^{\frac{4}{6}} +2t^{\frac{5}{6}}+8t^{\frac{6}{6}}+2t^{\frac{7}{6}}+2t^{\frac{8}{6}}+2t^{\frac{9}{6}}+t^{\frac{10}{6}}.
\end{equation} 
This  formula for the spectrum clearly implies
the following formula for the Alexander polynomial $\Delta^3(\A)$:
\begin{equation} 
\label{AlexNF.3}
\Delta^3(\A)=\Phi_1^{8}\cdot \Phi_2^{2}\cdot \Phi_3^{2}\cdot \Phi_6^{2}.
\end{equation} 
A generic plane section of $\A$ is a nodal line arrangement, and hence in this case $\Delta^1(\A)=\Phi_1^{5}$. It is easy to compute $\chi(M(\A))= -2$, and hence using
 the formula \eqref{Euler}, we get 
$$\Delta^2(\A)=\Phi_1^{10}.$$
Note that the two points $A=(0:0:0:1)$ and $B=(1:0:0:0)$ both correspond to dense edges $X$ with $n_X=4$, and any hyperplane in $\A$ contains at least one of these two points. It follows that
the value $k= 3$ is resonant with respect to the arrangement $\A$, i.e. the defining property of non-resonancy in Remark \ref{rkprim1} is not satisfied, but still there is no contribution to $H^m(F,\C)_{-1}$ for $m<3$.
This fact suggests that Yoshinaga's results in \cite{Yo0} for real line arrangements might have a higher dimensional analogue.
\end{ex}

\begin{ex}[A free discriminant surface]
\label{exD3}
Consider the surface in $\PP^3$ given by
$$V:f=y^2z^2-4xz^3-4y^3w+18xyzw-27x^2w^2=0.$$
Then $V$ is just the discriminant of cubic binary forms in $\PP(\C[u,v]_3)=\PP^3$, i.e. the set of cubic forms in $u,v$ with a multiple linear factor. It is known that $V$ is a free surface with exponents $d_1=d_2=d_3=1$ and $V$ is homeomorphic to $\PP^1 \times \PP^1$, see \cite{DStFS}.
Using the homogeneity under the obvious $G\ell_2(\C)$-action on $\PP(\C[u,v]_3)$, it is easy to see that $V$ is locally quasi-homogeneous.
Running the algorithm described above, we get the following for the terms occurring in the inequality  \eqref{eqHappy2}
\begin{equation} 
\label{spD3}
E_2^{3-q,q}(f)_k=0 \text{ for } (q,k) \ne (0,4) \text{ and  } \dim E_2^{3,0}(f)_4 =1.
\end{equation} 
It follows that $b_3(F) \leq 1$. On the other hand we have
$$\chi(F)=4 \chi(M)=4(\chi(\PP^3)-\chi(V))=4(4-4)=0.$$
Note that a generic plane section of $V$ is a quartic curve with 4 cusps, and hence $b_0(F)=1$ and $b_1(F)=0$, see for instance \cite[Proposition 4.4.8]{D1}. It follows that
$$b_3(F)=b_2(F)+b_0(F) \geq 1.$$
It follows that $b_3(F)=1$ and hence $f$ is $k$-top-computable for any $k\in [1,4]$. The corresponding Alexander polynomials are
$\Delta^3(V)=\Delta^0(V)=\Phi_1$ and $\Delta^2(V)=\Delta^1(V)=1$.
\end{ex}

\begin{ex}[A non-free irreducible surface]
\label{exGEN}
Consider the surface in $\PP^3$ given by
$$V:f=x^3z+x^2y^2+y^2w(y+w)=0,$$
and the corresponding Milnor fiber $F:f=1$ in $\C^4$. Let $H$ be the hyperplane in $\C^4$ given by $x=0$ and note that
$F_0=F \cap H$ is given by $y^2w(y+w)=1$ in $\C^3$, with coordinates $y,z,w$.
It follows that $F_0$ is a smooth surface, homotopically equivalent to the affine curve
$F_0':y^2w(y+w)=1$ in $\C^2$, with coordinates $y,w$. It is easy to see that the projective closure $C$ of $F_0'$ is a quartic irreducible curve, with a unique singular point of type $A_3$. It follows that $\chi(C)=2-(4-1)(4-2)+\mu(C)=-1$, see for instance \cite[Corollary 5.4.4]{D1}. Then $\chi(F_0')=\chi(C)-3=-4$,
since $C$ has 3 points at infinity. It follows that $b_0(F_0)=1$, $b_1(F_0)=5$ and $b_j(F_0)=0$ for $j \geq 2$. On the other hand, the projection on the $x$-coordinate induces a locally trivial fibration
$$ \C^2 \to F \setminus F_0 \to \C^*,$$
and hence $F \setminus F_0$ is homotopy equivalent to $\C^*$.
The Gysin sequence in homology, 
$$\cdots \to H_k(F\setminus F_0) \to  H_k(F) \to H_{k-2}(F_0) \to H_{k-1}(F\setminus F_0) \to \cdots,$$ see for instance \cite[Equation (2.2.13)]{D1}, yields the following Betti numbers for $F$: $b_0(F)=1$, $b_1(F)=b_2(F)=0$ and $b_3(F)=5$. Indeed, note that a generic plane section $V_P=V \cap P$ of $V$ is an irreducible curve in $P=\PP^2$ having a point of multiplicity $3$. It follows by \cite[Corollary 4.3.8]{D1} that 
$$\pi_1(P \setminus V_P)=\pi_1(\PP^3 \setminus V)=\Z/ d\Z.$$
Then \cite[Corollary 4.1.10]{D1} implies that $b_1(F)=0$.

Running the algorithm described above,  we get the following non-zero terms among those occurring in the inequality  \eqref{eqHappy2}
\begin{equation} 
\label{ssGEN}
\dim E_2^{2,1}(f)_k=1 \text{ for } k=1,2,3,4 \text{ and  } \dim E_2^{3,0}(f)_4 =1.
\end{equation}
It follows that we get an equality in \eqref{eqHappy2}, with implies that $V:f=0$ is $(k,1)$-top-computable for all the integers $k$ with $1 \leq k \leq 4$. In particular we have 
\begin{equation} 
\label{spGEN}
Sp^0_P(f)=t^{\frac{4}{4}} +t^{\frac{5}{4}}+t^{\frac{6}{4}}+t^{\frac{7}{4}}+t^{\frac{8}{4}}.
\end{equation} 
The corresponding Alexander polynomials are
$\Delta^3(V)= \Phi_1^2\cdot \Phi_2\cdot \Phi_4$, $\Delta^0(V)=\Phi_1$ and $\Delta^2(V)=\Delta^1(V)=1$. 

The value of the Alexander polynomial $\Delta^3(V)$ can also be obtained from the isomorphism
 $H_3(F) \to H_1(F_0)$ in the above Gysin sequence, using its functoriality, but not the spectrum $Sp^0_P(f)$. Note also that the  surface in $\PP^3$ given by
$$V':f'=x^3z+x^2y^2+y^2w(y+w)+x^2w^2=f+x^2w^2=0,$$
has the same topological properties as $V:f=0$, but this time the inequality in \eqref{eqHappy2}
is strict.
\end{ex}

\begin{ex}[Another non-free irreducible surface]
\label{exGEN2}
Consider the surface in $\PP^3$ given by
$$V:f=x^4+x^3z+yz^2w=0,$$
and the corresponding Milnor fiber $F:f=1$ in $\C^4$. Consider the hyperplanes $H_y:y=0$ and $H_z:z=0$ in $\C^4$ and set
$F_y=F \cap H_y$, $F_z=F \cap H_z$, $F_0=F_y \cap F_z$ and $F'=F\setminus F_0$.
It follows that $F_0$ is a smooth curve in $F$, homotopically equivalent to 4 points. 
The corresponding Gysin sequence in homology contains the sequence
$$0= H_4(F) \to H_{0}(F_0) \to H_{3}(F') \to H_3(F) \to H_{-1}(F_0)=0,$$
which implies $b_3(F)=b_3(F')-4$. Next $D'=(F_y \cup F_z) \setminus F_0$ is a smooth divisor in $F'$ and the projection
$$p: F' \setminus D' \to (\C^*)^2,  \  \   \ (x,y,z,w) \mapsto (y,z)$$
is a locally trivial fibration with contractible fibers and hence it is a homotopy equivalence.
Moreover, it is easy to see that $D'$ has five connected components, one homotopy equivalent to $\C$ minus 5 points, the other four homotopy equivalent to $\C^*$.
It follows that $b_1(D')=9$. The corresponding Gysin sequence of the pair $(F',D')$ contains the sequence
$$0=  H_{3}(F' \setminus D') \to H_3(F') \to H_{1}(D') \to H_{3}(F' \setminus D')=\Z,$$
which implies that $b_3(F') \in \{8,9\}$. Therefore we get $b_3(F) \in \{4,5\}$.

Running the algorithm described above,  we get the following non-zero terms among those occurring in the inequality  \eqref{eqHappy2}
\begin{equation} 
\label{ssGEN2}
\dim E_2^{2,1}(f)_k=1 \text{ for } k=1,2,3 \text{ and  } \dim E_2^{3,0}(f)_4 =1.
\end{equation}
The inequality  \eqref{eqHappy2} implies $b_3(F)=4$ and hence  $V:f=0$ is $(k,1)$-top-computable for all the integers $k$ with $1 \leq k \leq 4$. In particular we have 
\begin{equation} 
\label{spGEN2}
Sp^0_P(f)=t^{\frac{4}{4}} +t^{\frac{5}{4}}+t^{\frac{6}{4}}+t^{\frac{7}{4}}.
\end{equation} 
Using  \cite[Proposition 4.4.3]{D1}, we get as above $b_1(F)=0$. Consider the partition of $V$ given by $V=V_z \cup (V \setminus V_z)$, where $V_z=V \cap \{ z=0\}$.
Then $V_z=\PP^1$ and hence $\chi(V_z)=2$. Using the theory of tame polynomials, see for instance \cite{Br}, it follows that $\chi(V \setminus V_z)=2$, and hence 
$\chi(V)=\chi(V_z)+\chi(V \setminus V_z)=4$. It follows that $\chi(M)=\chi(\PP^3)-\chi(V)=0$.
This implies that corresponding Alexander polynomials are
$\Delta^3(V)= \Phi_1\cdot \Phi_2\cdot \Phi_4$, $\Delta^0(V)=\Phi_1$ and $\Delta^2(V)=\Phi_2\cdot \Phi_4$. In particular $b_2(F)=3$.

\end{ex}


\begin{thebibliography}{00}



\bibitem{A2}  E. Artal Bartolo, Topology of arrangements and position of singularities. Ann. Fac. Sci. Toulouse Math. (6), 23(2014), 223-265.








\bibitem{PB} P. Bailet,
On the monodromy of Milnor fibers of hyperplane arrangements,  Canad.\ Math.\ Bull.\ 57 (2014), 697--707.

\bibitem{BSett} P. Bailet, S. Settepanella, Homology graph of real arrangements and monodromy of Milnor fiber, arXiv:1606.03564.

\bibitem{BY} P. Bailet, M. Yoshinaga, Degeneration of Orlik-Solomon algebras and Milnor fibers of complex line arrangements, Geom.\ Dedicata 175 (2015), 49--56.

\bibitem{Br} S.A. Broughton,  Milnor numbers and the topology of polynomial
hypersurfaces. Invent. Math.{92} (1988), 217--241.

\bibitem{BS0} N. Budur and M. Saito,  Multiplier ideals, $V$-filtration, and spectrum, J. Alg. Geom. 14 (2005), 269--282.


\bibitem{BS} N. Budur and M. Saito, Jumping coefficients and spectrum of a hyperplane arrangement, Math. Ann.  { 347} (2010), 545--579.



\bibitem{BDS}   N. Budur, A. Dimca and M. Saito, First Milnor cohomology of hyperplane arrangements, Contemporary Mathematics  { 538} (2011), 279-292. 

\bibitem{Cal1} F. Callegaro, The homology of the Milnor fiber for classical braid groups. Algebr. Geom. Topol. 6 (2006), 1903--1923.


\bibitem{Cal2}  F. Callegaro,  Salvetti complex, spectral sequences and cohomology of Artin groups. Ann. Fac. Sci. Toulouse Math. (6) 23 (2014), no. 2, 267--296.

\bibitem{CMN}
F.J. Castro-Jim\' enez, D. Mond, L. Narv\'aez-Macarro,  Cohomology of the complement of a free divisor. Trans. Amer. Math. Soc. 348 (1996), 3037--3049. 




\bibitem{CS} D.~C.~Cohen, A.~I.~Suciu,
On Milnor fibrations of arrangements, 
J. London Math. Soc.  { 51} (1995), no.~2, 105--119.













\bibitem{Co} CoCoA-5 (15 Sept 2014): a system for doing Computations in Commutative Algebra, available at http://cocoa.dima.unige.it



\bibitem
{Sing} { W. Decker, G.-M. Greuel, G. Pfister \and H. Sch{\"o}nemann.} \newblock {\sc Singular} {4-0-1} --- {A} computer algebra system for polynomial computations, available at {http://www.singular.uni-kl.de} (2014).


\bibitem{De}  P.\ Deligne, Th\'eorie de Hodge II. Publ.\ Math.\ IHES { 40} (1972), 5--57.


\bibitem{Dcomp} A.~Dimca,   On the Milnor fibrations of weighted homogeneous polynomials, Compositio 
        Math. 76 (1990), 19-47.


\bibitem{D1} A.~Dimca, 
{\em Singularities and topology of hypersurfaces},
Universitext, Springer-Verlag, New York, 1992. 

\bibitem{D2} A. Dimca, {\em Sheaves in Topology},  Universitext, Springer-Verlag, 2004.

\bibitem{DNag} A.\ Dimca, Tate properties, polynomial-count varieties, and monodromy of hyperplane arrangements, Nagoya Math.\ J.\ 206 (2012), 75--97.




\bibitem{DAG}  A. Dimca, Jacobian syzygies, stable reflexive sheaves, and Torelli properties for projective hypersurfaces with isolated singularities, Algebraic Geometry 4(2017), 290--303.

\bibitem{DHA}  A. Dimca,   {\em Hyperplane Arrangements: An Introduction}, Universitext, Springer-Verlag, 2017


\bibitem{DL3} A.\ Dimca, G.\ Lehrer, Cohomology of the Milnor fiber of a hyperplane arrangement with symmetry, in: F.\ Callegaro et al.\ (eds.), Configuration Spaces, Springer INdAM Series 14, 2016, pp.\ 233--274.



\bibitem{DSVan} { A. Dimca,  M. Saito}, Some consequences of perversity of vanishing cycles, Ann.Inst. Fourier, Grenoble 54(2004), 1769--1792.


\bibitem{DS1} { A. Dimca,  M. Saito}, {  Koszul complexes and spectra of projective hypersurfaces with isolated singularities}, arXiv:1212.1081.

\bibitem{DS2} { A. Dimca,  M. Saito},  Generalization of theorems of Griffiths
and Steenbrink to hypersurfaces
with ordinary double points, arXiv:1403.4563, to appear in Bull. Math. Soc. Sci. Math. Roumanie.







\bibitem{DStEdin} A. Dimca, G. Sticlaru, { Koszul complexes and pole order filtrations, } Proc. Edinburg. Math. Soc. 58(2015), 333--354.



\bibitem{DStFor} A. Dimca, G. Sticlaru, A computational approach to Milnor fiber cohomology, Forum Math. 29 (2017),  831--846.

\bibitem{DStFS} A. Dimca, G. Sticlaru, Free and nearly free surfaces in $P^3$, arXiv:1507.03450. to appear in Asian J. Math.




\bibitem{DStMFgen} A. Dimca, G. Sticlaru, Computing the monodromy and pole order filtration on Milnor fiber cohomology of plane curves, arXiv: 1609.06818.
  



\bibitem{Lib} A.\ Libgober,  Eigenvalues for the monodromy
of the Milnor fibers of arrangements. In: Libgober, A., Tib\u ar, M.\ (eds) Trends in Mathematics: Trends in Singularities. Birkh\" auser, Basel (2002)

\bibitem{MSS} L. Maxim, M.Saito and J. Sch\"urmann, Spectral Hirzebruch-Milnor classes of singular hypersurfaces,
arXiv:1606.02218.


\bibitem{MP}
A. M\u acinic and S. Papadima, On the monodromy action on Milnor fibers of graphic arrangements, Topology Appl.  { 156} (2009), 761-774.


\bibitem{MPP}  A. M\u acinic, S. Papadima, C.R. Popescu, Modular equalities for complex reflexion arrangements, Documenta Math. 22 (2017) 135--150.


\bibitem{NM} L. Narv\' ez Macarro, A duality approach to the symmetry of Bernstein-Sato polynomials for free divisors, Adv. Math. 281 (2015), 1242--1273.





\bibitem{OT}P.~Orlik and H.~Terao, Arrangements of Hyperplanes,
Springer-Verlag, Berlin Heidelberg New York, 1992.



\bibitem{PS}  S. Papadima, A. I. Suciu, The Milnor fibration of a
hyperplane arrangement: from modular resonance to algebraic monodromy,  
Proc.\ London Math.\ Soc.\ DOI: 10.1112/plms.12027.

\bibitem{PeSt} C.\ Peters, J.\ Steenbrink, {\em Mixed Hodge Structures},
Ergeb.\ der Math.\ und ihrer Grenz.\ 3. Folge 52,
Springer, 2008.



\bibitem{SaSurvey} M. Saito, On b-function, spectrum and multiplier ideals. Algebraic analysis and around, 355--379, Adv. Stud. Pure Math., 54, Math. Soc. Japan, Tokyo, 2009.





\bibitem{Sa1} M. Saito, Multiplier ideals, b-function, and spectrum of a hypersurface singularity, Compositio Math.143 (2007), 1050--1068.

\bibitem{Sa0} M.\ Saito, Bernstein--Sato polynomials of hyperplane arrangements, Selecta Math.\ {22} (2016),  2017--2057.

\bibitem{Sa2} M. Saito, Hilbert series of graded Milnor algebras and roots of Bernstein-Sato polynomials, arXiv: 1509.06288.

\bibitem{Sa3} M. Saito, Bernstein-Sato polynomials and graded Milnor algebras for projective hypersurfaces with weighted homogeneous isolated singularities, arXiv:1609.04801.

\bibitem{Sa4} M. Saito, Roots of Bernstein-Sato polynomials of homogeneous polynomials with 1-dimensional  singular loci, arXiv:1703.05741.

\bibitem{SalSet} M. Salvetti, S. Settepanella, 
Combinatorial Morse theory and minimality of hyperplane arrangements, 
Geom. Topol. 11 (2007), 1733--1766. 


\bibitem{Se2} S. Settepanella,  Cohomology of pure braid groups of exceptional cases. Topology Appl. 156 (2009), no. 5, 1008--1012.


\bibitem{Steen} J. Steenbrink, Intersection form for quasi-homogeneous singularities, Compositio Math. 34(1977), 211-223.



\bibitem{S2} A. Suciu,
Hyperplane arrangements and Milnor fibrations, Annales de la Facult\'e des Sciences de Toulouse 23 (2014), no. 2, 417-481.

\bibitem{Tapp} D. Tapp, Picard-Lefschetz monodromy of products, J. Pure Appl. Algebra 212 (2008),
2314--2319.


\bibitem{WY}  J. Wiens, S. Yuzvinsky,  De Rham cohomology of logarithmic forms on arrangements of hyperplanes. Trans. Amer. Math. Soc. 349 (1997), 1653--1662.

\bibitem{Terao78} H.\ Terao, Forms with logarithmic pole and the filtration by the order of the pole, Proc.\ Intern.\ Symposium Algebraic Geometry, Kyoto 1977, pp.\ 673--685, Kinokuniya Book Store, Tokyo, 1978.

\bibitem{Te} H. Terao, Generalized exponents of a free arrangement of hyperplanes and Shephard-Todd-Brieskorn formula, Invent. math. 63(1981), 159-179.





\bibitem{Yo0} M. Yoshinaga, Milnor fibers of real line arrangements, Journal of Singularities, vol 7 (2013), 220-237.


\bibitem{Yo} M. Yoshinaga, Freeness of hyperplane arrangements and related topics, Annales de la Facult\'e des Sciences de Toulouse, vol. 23 no. 2 (2014), 483-512.























\end{thebibliography}
\end{document}